\documentclass[12pt]{article}
\usepackage{amsmath,amsthm,amsfonts,amssymb,amscd}
\def\qed{{\unskip\nobreak\hfil\penalty50
\hskip2em\hbox{}\nobreak\hfil$\square$
\parfillskip=0pt \finalhyphendemerits=0\par}\medskip}
\def\proof{\trivlist \item[\hskip \labelsep{\bf Proof\ }]}
\def\endproof{\null\hfill\qed\endtrivlist}
\def\End{{\mathrm {End}}}
\def\Hom{{\mathrm {Hom}}}

\def\dim{{\mathrm {dim}}}

\def\e{\varepsilon}

\def\phi{\varphi}

\def\Om{\Omega}

\def\r{{\rho}}

\def\emptyset{\varnothing}
\def\setminus{\smallsetminus}

\def\Diff{{\mathrm {Diff}}}
\def\Mob{{\rm\textsf{M\"ob}}}

\def\End{{\mathrm {End}}}
\def\Hom{{\mathrm {Hom}}}

\def\dim{{\mathrm {dim}}}

\def\e{\varepsilon}

\def\phi{\varphi}

\def\Om{\Omega}

\def\r{{\rho}}

\newtheorem{theorem}{Theorem}[section]
\newtheorem{lemma}[theorem]{Lemma}
\newtheorem{conjecture}[theorem]{Conjecture}
\newtheorem{corollary}[theorem]{Corollary}
\newtheorem{definition}[theorem]{Definition}

\newtheorem{proposition}[theorem]{Proposition}
\newtheorem{remark}[theorem]{Remark}

\def\emptyset{\varnothing}
\def\setminus{\smallsetminus}

\def\Ind{{\mathrm {Ind}}}
\def\Diff{{\mathrm {Diff}}}
\def\Mob{{\rm\textsf{M\"ob}}}

\def\A{{\cal A}}

\def\B{{\cal B}}

\def\I{{\cal I}}

\def\H{{\cal H}}

\renewcommand{\qed}{\ \hfill $\blacksquare$}

\newcommand{\bdefin}{\begin{definition}}
\newcommand{\blemma}{\begin{lemma}}
\newcommand{\bprop}{\begin{proposition}}
\newcommand{\btheor}{\begin{theorem}}
\newcommand{\bcoro}{\begin{corollary}}
\newcommand{\bconj}{\begin{conjecture}}
\newcommand{\edefin}{\end{definition}}
\newcommand{\elemma}{\end{lemma}}
\newcommand{\eprop}{\end{proposition}}
\newcommand{\etheor}{\end{theorem}}
\newcommand{\ecoro}{\end{corollary}}
\newcommand{\econj}{\end{conjecture}}
\newcommand{\brem}{\begin{remark}}
\newcommand{\erem}{\end{remark}}

\newcommand{\ba}{\begin{array}}
\newcommand{\ea}{\end{array}}
\newcommand{\bea}{\begin{eqnarray}}
\newcommand{\eea}{\end{eqnarray}}
\newcommand{\bean}{\begin{eqnarray*}}
\newcommand{\eean}{\end{eqnarray*}}

\title{\huge Intersections of finite families of finite index subfactors.\\}
\author{
{\sc Vaughan F.R. Jones}\footnote{Supported in part by NSF gramt DMS-9322675 and Marsden
grant UOA520.}\\
Department of Mathematics\\
U.C.Berkeley\\
Berkeley, CA 94720\\
E-mail: {\tt vfr@math.berkeley.edu}\\
{}\\
{\sc Feng Xu}\footnote{Supported in part by NSF grant DMS-0200770.}\\
Department of Mathematics\\
University of California at Riverside\\
Riverside, CA 92521\\
E-mail: {\tt xufeng@math.ucr.edu}}
\begin{document}
\date{}
\maketitle

\begin{abstract}
We prove that finiteness of the 
index of the intersection of a finite set of finite index subalgebras 
in a von Neumann algebra (with small centre) is equivalent to the
finite dimensionality of the algebra generated by
the conditional expectations onto the subalgebras. \par
\end{abstract}
%\thanks
%2000 Mathematics Subject Classification. 46S99, 81R10. 
%Key words: Subfactors, Conformal field theories.
%\endthanks 
\newpage

\section{Introduction}
Subfactor theory provides an entry point into a world of mathematics and
physics containing large parts of conformal field theory, 2-dimensional statistical mechanical
models, quantum groups, finitely presented groups and low dimensional topology.
We cite two facts attesting to the effectiveness of the subfactor point of view.
The first is the discovery of a polynomial invariant of links in \cite{J3}
and the second is the use of subfactors by Wassermann in \cite{W} to give a
rigorous and entirely unitary definition of fusion of loop group representations.

But the very breadth of areas in this world means that one would not expect
subfactors to yield the most effective technique for
any of the particular topics mentioned above. As things stand for instance it
would be foolish to expect subfactor ideas to have a big impact on 
finite group theory. It is more reasonable to expect subfactors to provide the
machinery for theories close to the abstract structure of a subfactor itself. 
The most obvious such structure is an intermediate subfactor and the 
first author and Bisch in \cite{BJ} led to the discovery of an entirely new 
algebra, in the spirit of the type A Hecke algebra,  which they called the
"Fuss-Catalan" algebra. This success is a compelling argument for further work 
on the lattice of intermediate subfactors of a factor, which had been suggested in
papers by Watatani (\cite{Wat})  and Watatani and Sano (\cite{SW}) . The main contribution of
\cite{SW} was to introduce the angle operator between two subfactors as an invariant
of a pair of subfactors and calculate several examples. The following question 
is hinted at but never mentioned explicitly in \cite{SW}:

$\cal Q$:{\it If the spectrum of the angle operator between two finite index subfactors is finite,
is the intersection of finite index? }

(Note that it follows from the elementary theory of the index that the angle condition
is necessary for the intersection to have finite index.) The simplest non-trivial example
of a subfactor is given by a pair $\Gamma_0 \subseteq \Gamma$ of discrete groups each having
infinite non-identity conjugacy  classes. The group von Neumann algebras give the factor 
and subfactor $M$ and $N$ and $[M:N]=[\Gamma:\Gamma_0]$. Of course the intersection of
two finite index subgroups has finite index so this case is rather special. Indeed what
is special is that the orthogonal projections (=conditional expectations) onto these subgroup
subfactors all commute so that angles are all either $0$ or $\pi/2$ so the angle operator has
only two elements in its spectrum. Another simple construction of subfactors is as fixed point
algebras. If $G$ and $H$ are two finite groups of automorphisms of a II$_1$ factor $M$ then
the fixed point algebras $M^G$ and $M^H$ are of finite index and the conditional expectations
will not commute so that $\cal Q$ becomes quite relevant. The simplest example is when 
$G$ and $H$ are are $\mathbb Z /2\mathbb Z$ so that the group generated by them is dihedral. 
It is easy to check that the answer to $\cal Q$ is affirmative in this case. This example,
and question $\cal Q$, make perfect sense in Galois theory: if $E$ and $F$ are two subfields
of a field $K$ with $[K:E]$ and $[K:F]$ finite there are trace maps $Tr_E$ and $Tr_F$ which
play the role of the conditional expectations and one may ask if  finite rank
of the $\mathbb Z$-algbera generated by $Tr_E$ and $Tr_F$ guarantees that 
$[K:E\cap F]<\infty$.
George Bergmann (\cite{Berg}) answered this question in the affirmative in zero characteristic,
thus reinforcing $\cal Q$. He further extended it to a finite collection of subfields rather
than just two. Similarly we can extend $\cal Q$ to $\cal Q$':

$\cal Q$' {\it If the conditional expectations onto a finite family of finite index subfactors
generate a finite dimensional algebra, is the intersection of the subfactors of finite index?}

Of course there is no reason for the intersection of subfactors to be a factor so 
we must use some alternative definition of finite index to properly formulate 
$\cal Q$' in detail. Several equivalent definitions are available, the most general
being the probabilistic one of Pimsner and Popa in \cite{PP}. In the case of finite von Neumann
algebras one may use the property that the large algebra is a finitely generated
left (or right) module over the small one. For properly infinite  algebras we will use
the endomorphism theory pioneered by Longo (see \cite{L1}). In this paper  we will answer
$\cal Q$' in the affirmative. Our proof in the properly infinite case
uses the properties of type III factors and has applications to conformal
field theories where type III factors appear naturally.\par

The paper is organized as follows: After introducing some
basics of index theory and setting up notations in \S2, in \S3 we prove
$\cal Q$' in   Th. \ref{finitecase}, and we give two 
applications in Cor. \ref{3.2} and Th. \ref{3.3}. In \S4 we first prove 
Th. \ref{key} using two lemmas. An averaging technique in \cite{ILP}
plays a key role in the proof of Lemma \ref{s=1}. Th. \ref{key} implies
Cor. \ref{general} which proves the extension of $\cal Q$'. In \S4.3,
we describe the setting of conformal nets where the assumptions of
Cor. \ref{general} are naturally satisfied (cf. Lemma \ref{irred1}, 
Cor. \ref{subnet1} and Lemma \ref{irred2}), and apply Cor. \ref{general} to
a large class of conformal nets in Cor. \ref{subnet2}.

The first author would like to thank D. Bisch for useful discussions and motivation
for the problem.
   
\section{Preliminaries}
Index theory of subfactors was initiated in \cite{J1} in 
the setting of finite factors. In this paper we will use the following
more general definition from \cite{PP} and \cite{P}:
\begin{definition}\label{ind}
Let $N\subset M$ be an inclusion of von Neumann algebras with a 
conditional expectation $E:M\rightarrow N$. The (probabilistic) index
of $E$, denoted simply by $\Ind E$, is defined by
$$
\Ind E= (sup \{c\geq 0| E(m)\geq cm, \forall m\in M_+\})^{-1}.
$$ 
The inequality in the above definition will be referred to as 
Pimsner-Popa inequality. 
The inclusion  $N\subset M$ has finite index if there exists a 
conditional expectation $E:M\rightarrow N$ such that 
$\Ind E<\infty$.
\end{definition}
Let $N\subset M$ be an inclusion of von Neumann algebras and $E:
M\rightarrow N$ a normal faithful conditional expectation. Let $\phi$
be a normal faithful state on $M$ such that $\phi=\phi\cdot E$. Let
$\H$ be the Hilbert space of the GNS representation associated to
$\phi$, and let $\Omega\in \H$ be the vector such that
$\phi(m)=\langle m\Omega, \Omega
\rangle, \forall m\in M.$  Let $e$ be the Jones projection from 
$\H$ to $\overline{N\Omega}$.  
\begin{definition}\label{basis}
A family of elements $\{ m_j\} \subset M$ satisfying the conditions:\par
(1) $E(m_i^*m_j)= \delta_{ij} f_j$ where $f_j$ is a projection in $N, 
\forall i,j$;\par
(2) $\overline{\sum_j m_j(e\H)}=\H$ \par
is called an orthonormal basis of $M$ over $N$ via $E$.
\end{definition} 

\begin{lemma}\label{basic}
Let $M\subset B(\H)$ be a von Neumann algebra represented standardly 
on a Hilbert space $\H$ with a cyclic separating vector $\Omega$ and
denote by $\phi$ the vector state on $M$ with 
$\phi(m)=\langle m\Omega, \Omega
\rangle.$
Let $N\subset M$ be a von Neumann algebra and assume that $E:M\rightarrow N$
is a faithful conditional expectation with $\phi(E(m))= \phi(m), 
\forall m\in M.$ Let $e$ be the Jones projection from $\H$ to 
$\overline{N\Omega}$, and $J$ the canonical conjugation of $M$ 
with respect to $\Omega$. Denote by $\langle M, e\rangle$ the von Neumann
algebra generated by $M, e$.
Then: \par
(1) $MeM$ is weakly dense in $\langle M, e\rangle$; \par
(2) $\langle M, e\rangle= JN'J$; \par
(3) $N\subset M$ has finite index if and only if $M\subset \langle M, e\rangle$
has finite index; \par
(4) $N\subset M$ has finite index if and only if 
$N\otimes B \subset M\otimes B$ has finite index where $B$ is a factor.
\end{lemma}  
\proof
(1), (2) is contained in \cite{J1}. (3) follows by 1.2 of Page 9 of
\cite{P}. As for (4), note that if $E: M\otimes B\rightarrow N\otimes B$ is
a conditional expectation which verifies Pimsner-Popa inequality for
some constant $c>0$, then $E(M\otimes 1)= (N\otimes1)$, hence by 
restriction $N\subset M$ has finite index. On the other hand if
$N\subset M$ has finite index and assume that 
$E: M\rightarrow N$ is
a conditional expectation which verifies Pimsner-Popa inequality for
some constant $c'>0$, by Th. 1.1.6 of \cite{P} there is an orthonormal 
basis $\{ m_i\}$ of $M$ over $N$ via $E$, such that 
$\sum_jm_jm_j^*$ is bounded.  By definition  
$\{ m_i\otimes 1\}$  is
also an orthonormal basis of  $M\otimes B$ over $N\otimes B$ via 
$E\otimes id$, and hence   by Th. 1.1.6 of \cite{P} again 
$N\otimes B \subset M\otimes B$ has finite index.   
\endproof

\section{Finite type von Neumann algebra case}\label{finite}
Let $M$ be a direct sum of finitely many finite
factors with faithful trace $tr$.  Let $\cal P$ be a finite set of (unital)
finite index subalgebras of $M$ and for each $P\in \cal P$
let $e_P$ be the projection from $L^2(M,tr)$ onto $L^2(P,tr)$.
Let ${\cal F}= \{e_P : P\in {\cal P}\}$.

\begin{theorem} \label{finitecase}
If $N=\bigcap_{P\in {\cal P}} P$ then 

$$[M:N]<\infty \iff \dim \{{\cal F}''\} <\infty.$$
\end{theorem}

\proof (only if) If $[M:N]<\infty$ then $\dim \{ N'\cap \langle M,e_N \rangle\}<\infty$ 
 by reduction to the factor case. But ${\cal F}'' \subseteq  N'\cap \langle M,e_N \rangle$.
\vskip 4pt
(if)  The key thing to establish is that $Me_N M$ is contained in a finitely 
generated
left $M$-module contained in $\langle M,e_N \rangle$.

Since ${\cal F}''$ is finite dimensional and $e_N =\inf_{f\in {\cal F}} f$,
$e_N$ is a polynomial in the $\{f\in \cal F\}$. Thus 
$Me_NM \subseteq \sum_{w\in {\cal W}} MwM$ where $\cal W$ is a finite set of
words on $\{f\in \cal F\}$. Now let $w \in {\cal W}$ be $f_1f_2...f_k$.
Then for each $i$, $Mf_iM$ (which is contained in $\langle M,e_N \rangle$)
is a basic construction for a finite index
subalgebra so by [PP],[GHJ] it is a finitely generated left $M$-module.
Hence $Mf_1Mf_2Mf_3...f_kM$ is also a finitely generated left $M$-module,
being the image under the multiplication map of 
$(Mf_1M) \otimes_M (Mf_2M)\otimes_M (Mf_3M)\otimes_M... \otimes_M (Mf_kM)$.
Thus since $MwM\subseteq Mf_1Mf_2Mf_3...f_kM$, $Me_NM$ is contained in
a finitely generated left $M$-module inside $\langle M,e_N \rangle$, say $$Me_NM\subseteq 
\sum
_{i=1}^n Mr_i$$
for $r_i \in \langle M,e_N \rangle$.

Now consider $\langle M,e_N \rangle$. It has a faithful normal trace $Tr$ 
so we may consider its action on $L^2(\langle M,e_N \rangle, Tr)$. By contradiction 
we suppose $\langle M,e_N \rangle$ is not a finite von Neumann algebra. 

First suppose $M_*$ is separable,
then so is $L^2(\langle M,e_N \rangle, Tr)$ so by proposition 3.14 of \cite{Tak}
there is a cyclic vector for  $Me_NM$ in $L^2(\langle M,e_N \rangle, Tr)$, call it $\Omega$.
Then $\{r_i \Omega\}$ is a finite set which is cyclic for $M$.
By the coupling constant, $M'$ is finite (on $L^2(\langle M,e_N \rangle, Tr)$) 
and $\langle M,e_N \rangle\subseteq JM'J$.

The nonseparable case requires a little more care. We know that $\langle M,e_N \rangle$ is 
a semifinite von Neumann algebra whose center, being the same as the center 
of $N$, has only countably many mutually orthogonal projections. Thus we may,
by   proposition 1.40 of \cite{Tak}, find an infinite, $\sigma$-finite 
projection
$p\in \langle M,e_N \rangle$ of central support $1$ with $p=\sum_n p_n$ where $p_n$ are
mutually orthogonal projections with $Tr(p_n)<\infty$. Choose positive numbers
$a_n$ with $\sum_n a_n^2 Tr(p_n) <\infty$. 
Then $\Omega=\sum_n a_n p_n$ is a vector in $L^2(\langle M,e_N \rangle,Tr)$
and it is a separating vector for $p\langle M,e_N \rangle p$ on $pL^2(\langle M,e_N \rangle,Tr)$. 
Since the central support of $p$ is $1$, the action of $\langle M,e_N \rangle$ on 
$L^2(\langle M,e_N \rangle,Tr)p$
is faithful so as in the previous paragraph  $\{r_i \Omega\}$ is a finite 
set which is cyclic for $M$.
Hence by the coupling constant, $M'$ is finite (on $L^2(\langle M,e_N \rangle,
 Tr)p$) 
and it contains the infinite von Neumann algebra 
$Jp\langle M,e_N \rangle pJ\quad -\quad$a contradiction. (Note that
$J$ is the canonical involution on $L^2(\langle M,e_N \rangle, Tr)$.)
\endproof

If we know that $N$ and $M$ are II$_1$ factors with $[M:N] <\infty$, the 
proof of \ref{finitecase}
gives a bound on $[M:N]$ in terms of $\dim A$ 
(where we set $A={\cal F}''$) and 
the individual
$[M:P]$ for $P\in \cal P$. To give an explicit (but rather crude) bound, let 
$\displaystyle L=\max_{P\in {\cal P}} [M:P]$ and $\ell$ be the length of 
the longest 
word $w$ in some basis of $A$ consisting of words on $\cal F$.

\begin{corollary}\label{3.2}
With hypotheses as above, $$[M:N] \leq L^\ell \hspace{3pt} \dim A.$$

\end{corollary}

\proof Since $\langle M,e_N \rangle \quad \subseteq \quad \sum_w MwM$ it suffices to bound 
$\dim_M(\overline{MwM})$ for each $w$. But since all indices are finite one 
may 
take the Connes
tensor product of the $L^2(MfM)$ instead of the algebraic one and
 $\dim_M(\overline{MwM}) \leq L^{length(w)}$ follows from the 
multiplicativity of the 
 $M$-dimension under tensor product and the fact that $MfM$ is just a 
basic construction for 
 a subactor in $\cal P$,so its $M$-dimension is bounded by $L$. 
\endproof 
Using the formalism of planar algebras (cf.\cite{J2}) the first author has 
found the slightly better bound 
(when $|{\cal P}|=2$),
$$[M:N] \leq \frac{L^\ell \hspace{3pt} \dim A}{4}.$$ 
In practice the bound seems to be
a lot sharper. \par
Given a subfactor it is always possible to perturb it 
(for instance conjugating it
by a unitary) to obtain a pair of subfactors to which \ref{finitecase} can be
applied. We have only carried this analysis out in one case, 
namely the locally
trivial subfactors coming from automorphisms which can be thought of as 
perturbations
of the tensor product subfactor. We obtain the following result which does not
seem easy to prove by other means.

\begin{theorem} \label{3.3}
Let $\cal A$ be a finite set of automorphisms of the finite 
factor
$M$. Then $N^{\cal A} = \{x\in N | \alpha(x)=x \hspace{3pt}\forall \alpha 
\in {\cal A}\}$ is of
finite index in $N$ iff the spectrum of the operator
 $$T=\sum_{\alpha, \beta \in {\cal A}}\alpha \beta^{-1}$$
is finite. 
\end{theorem}
\proof 
Let $M$ be the finite direct sum of factors 
$\oplus_{\alpha \in {\cal A}} N$
 where we take a copy of $N$ for each $\alpha$, and 
let $P$ and $Q$ be the
 two subfactors $P=\{\oplus_\alpha x |x\in N\}$ and 
$Q=\{\oplus_\alpha \alpha(x)
|x\in N\}$.
 Then $P\cap Q =\{\oplus_\alpha x |x\in N^{\cal A}\}$ so that 
the index of $N^{\cal A}$ in
 $N$ is finite iff the index of $P\cap Q$ in $N$ is finite. 
Use the obvious trace on
 $M$ to form conditional expectations. Identify $P$ with $N$ in the 
obvious way. We have 
$\displaystyle E_P(\oplus_\alpha x_\alpha)= 
\frac{\sum_\alpha x_\alpha}{|{\cal A
}|}$
and 
$$ E_Q(\oplus_\alpha x_\alpha)= 
\oplus_{\beta \in {\cal A}} \quad
\beta(\frac{\sum_\alpha \alpha^{-1}(x_\alpha)}{|{\cal A}|}).$$
So we see that $E_PE_QE_P$ is, up to a multiple, our operator $T$. Now
the algebra generated by two idempotents $p$ and $q$ is finite dimensional iff
the algebra generated by $pqp$ is finite dimensional, which in this context
is the same as the finiteness of the spectrum of $pqp$. So by \ref{finitecase}
we are done.
\endproof

The result appears to be non-trivial especially if the 
automorphisms $\alpha$ are inner,
in which case it states that the algebra generated by 
$Ad (u)$ for some finite set
of unitaries $u$ is finite dimensional iff 
$\sum_{u,v} Ad (uv^*)$ has finite spectrum.
Note that it is NOT true that the algebra generated by $u$ is 
finite dimensional iff $\sum_{u,v}  uv^*$ has finite spectrum.

\section{General factor case (with separable predual)} 
\subsection{Preliminaries on Sectors}
Let $M$ be a properly infinite factor.  The {\it sectors of $M$}  are given by
$$\text{Sect}(M) = \text{End}(M)/\text{Inn}(M),$$
namely $\text{Sect}(M)$ is the quotient of the semigroup of the
endomorphisms of $M$ modulo the equivalence relation: $\rho,\rho'\in
\text{End}(M),\, \rho\thicksim\rho'$ iff there is a unitary $u\in M$
such that $\rho'(x)=u\rho(x)u^*$ for all $x\in M$.

$\text{Sect}(M)$ is a $^*$-semiring (there is an addition, a product and
an involution) equivalent to the Connes correspondences (bimodules) on
$M$ up to unitary equivalence. If $\r$ is
an element of $\text{End}(M)$ we shall denote by $[\r]$
its class in $\text{Sect}(M)$. The operations are:

{\it Addition} (direct sum): Let $\rho_1,\rho_2,\dots\r_n\in\text{End}(M)$.
Choose a non-degenerate $n$-dimensional Hilbert $H$ space of isometries in $M$
and a basis $v_1,\dots v_n$ for $H$. Here $H$ is non-degenerate means that
the isometries $v_1,\dots v_n$ verifies $\sum_i v_iv_i^*= 1$. 

Then
\[
\r(x)\equiv\sum_{i=1}^n v_i\r_i(x)v^*_i,\quad x\in M,
\]
is an endomorphism of $M$. The definition of the direct sum 
endomorphism $\r$ does not depend on the choice of $H$ or on the 
basis, up to inner automorphism of $M$,
namely $\r$ is a well-defined sector of $M$.

{\it Composition} (monoidal product). The usual composition of maps
$$\rho_1\cdot\rho_2(x) = \rho_1(\rho_2(x)), \qquad x\in M,$$ defined on
$\text{End}(M)$ passes to the quotient $\text{Sect}(M)$.
Let $\rho \in \text{\rm End}(M)$ and $\e$ be a normal
faithful conditional expectation
$\e:
M\rightarrow \rho(M)$.  We define a number $d_\e\geq 1$ (possibly
$\infty$) by:
$$
d_\e^{-2} :=\text{\rm Max} \{t\in [0, +\infty)|
\e (m_+) \geq t m_+, \forall m_+ \in M_+
\}$$ (Pimsner-Popa inequality in \cite{PP}).\par
 We define
$$
d(\rho) = \text{\rm Min}_\e \{ d_\e \},
$$ 
where the minimum is taken over $\e$ with $d_\e < \infty$ (otherwise 
we put $d(\r) =\infty$).  $d(\rho)$ is called the  dimension of  $\rho$. 
We say that $\rho$ has finite index if $d(\rho)< \infty.$
It is 
clear from the definition that  the dimension  of  $\rho$ 
depends only the sector $[\rho]$.
The following properties of the dimension can be found in
\cite{L1}.
\begin{lemma}\label{dim}
(1) If $[\rho]=[\rho_1]+[\rho_2],$ then $ 
d(\rho)= d(\rho_1)+d(\rho_2);$\par
(2) $d(\rho_1\rho_2)= d(\rho_1) d(\rho_2).$
\end{lemma}
For $\lambda $, $\mu \in \text{\rm End}(M)$, we will use
$\text{\rm Hom}(\lambda , \mu )$ denote the vector space of intertwiners from
$\lambda $ to $\mu $, i.e. $a\in \text{\rm Hom}(\lambda , \mu )$ iff
$a\in M, a \lambda (x) = \mu (x) a $ for any $x \in M$.
A sector $\lambda$
is said to be irreducible if the vector space 
$\text{\rm Hom}(\lambda , \lambda )$ has
dimension one.\par
\begin{lemma}\label{can}
Let $N\subset M$ be a properly infinite subalgebra and assume that there
exists a normal faithful conditional expectation $E: M\rightarrow N$.
Suppose that M is represented standardly on a Hilbert space $H$ with 
a cyclic separating vector $\Omega,$  and the vector state 
$\omega (m)=\langle m\Omega,\Omega\rangle$ on $M$ verifies 
$\omega (E(m))=\omega (m), \forall m\in M.$ Let $e$ be the Jones projection
from $H$ onto $\overline{N\Omega}$ and $M_1:=\langle M,e\rangle$ be 
the von Neumann algebras generated by $M,e$. Then:\par
(1) Let $v\in M_1$ be an isometry with $vv^*=e.$ Then for every 
$m_1\in M_1$, there is a unique element, denoted by $\gamma(m_1) \in N,$
such that $vm_1 = \gamma(m_1) v$, and $\gamma \in  \text{\rm End}(M_1)$;\par
(2) We will use the same notation $\gamma$ to denote its restriction to 
$M$. Then $N\subset M$ has finite index if and only if 
$\gamma(M)\subset M$ has finite index.
\end{lemma}
\proof
The first part is Prop. 2.9 of \cite{LR}, and the second part follows from
(3) of Lemma \ref{basic}.
\endproof
\begin{remark}
The endomorphism $\gamma$ in (2) of Lemma \ref{can} is called 
canonical endomorphism for $N\subset M$ in \cite{LR}.
\end{remark} 
\subsection{General factor case}
\begin{theorem}\label{key}
Let $M$ be a factor with separable predual, and let 
$N\subset M$, $E:M\rightarrow N, \H , \Omega$ and $e$ 
be as in Lemma \ref{basic}.
If $e=\sum_{1\leq i\leq k} m_i R_i$ with $m_i\in M, R_i\in B(\H)$ such that
$R_i m = \rho_i(m) R_i, \forall m\in M$ and each $\rho_i\in  \text{\rm End}(M)
$ has finite index, $i=1,...,k.$ Then
$N\subset M$ has finite index.
\end{theorem}
\begin{remark}
Note that  each $\rho_i\in  \text{\rm End}(M)
$ has finite index, $i=1,...,k $ is a necessary condition when 
$N\subset M$ are properly infinite von Neumann algebras,  since
by Lemma \ref{can} $e= vv^*= \gamma(v^*) v, vm=\gamma(m)v , \forall m\in M$ 
for 
any pair $N\subset M$ as in Lemma  \ref{can}, including the case
when $N\subset M$ has infinite index. Also note that the conditions 
on Jones projection $e$ in Theorem \ref{key} are similar to that 
of Theorem \ref{finitecase}, but neither theorem implies the other.
\end{remark}
The proof of this theorem is divided into following steps consisting of
two lemmas. \par
First note that replacing 
$N,M, E, e, m_i, R_i, \rho_i$ by
$N\otimes B,M\otimes B, E\otimes id, e\otimes 1, 
m_i\otimes 1, R_i\otimes 1, \rho_i\otimes id$ respectively if necessary 
where $B$ is a type III factor with
separable dual, by (4) of Lemma \ref{basic} it is enough to prove
the theorem for $N\otimes B\subset M\otimes B$. So 
we can assume that $M$ is a type III factor, and $N$ is a type III
von Neumann algebra. Let $v, \gamma$ be as in Lemma \ref{basic} such
that $e=vv^*= \gamma(v^*)v, vm=\gamma(m)m, \forall m\in M.$  
Note that since each $\rho_i\in  \text{\rm End}(M)
$ has finite index, $i=1,...,k$, $\rho_i$ can be decomposed into 
sum of finitely many irreducible sectors. Hence we can assume that
$e=\sum_{1\leq i\leq n} \sum_{1\leq \alpha\leq l(i)}
m_{i\alpha} R_{i\alpha}$ with $m_{i\alpha}\in M, 
R_{i\alpha}\in B(\H)$ such that
$R_{i\alpha} m = \rho_i(m) R_{i\alpha}, \forall m\in M$ and each 
$\rho_i\in  \text{\rm End}(M)
$ is irreducible , has finite index, 
$1\leq \alpha\leq l(i) <\infty$ is a label, and
$[\rho_i]\neq [\rho_j], {\mathrm if} i\neq j, i,j=1,...,n.$ Moreover, 
replacing $ \{R_{i\alpha}\}$ by a maximal linearily independent 
subset  over ${\mathbb C}$ if necessary, 
we can
assume that $ \{R_{i\alpha}\}$ are linearily independent over ${\mathbb C}$.
\begin{definition}\label{V}
Let $V\subset M$ be a vector space consisting of all finite linear
combinations of 
elements $s_i$ with the property that
$s_i ^* v = \sum_{1\leq \alpha \leq l(i)} c_\alpha R_{i\alpha},
c_\alpha\in {\mathbb C}, i=1,...,n.$ Denote by 
$W$ the vector space spanned by $\{R_{i\alpha}\}$.
\end{definition}
\begin{lemma}\label{Hil}
$V$ is a finite dimensional Hilbert space with natural 
inner product $\langle s, t\rangle = t^* s$, and when $V$ is not
zero we can choose an orthonormal basis $\{s_{i\beta}\}$ with the property
$s_{i\beta} \in \Hom(\rho_i,\gamma), s_{i\beta}^* s_{i'\beta'}= 
\delta_{ii'}\delta_{\beta\beta'}.$
\end{lemma}
\proof
Let $W$ be the vector space spanned by $\{R_{i\alpha}\}$. Then 
$W$ is a vector space with dimension $\prod_{1\leq i\leq n} l(i)$. 
Let $F:V\rightarrow W$ be a conjugate linear map defined by
$F(s)= s^*v$. We claim that $F$ is one-to-one: if
$s^*v=0$, then $s^* vv^*= s^*e=0$, and so
$s^* \Omega=0$ which implies that $s=s^*=0$ since $\Omega$ is
separating for $M$. Hence $dim V\leq dim W$.\par
Assume that $s_i\in V, s_i ^* v = \sum_{1\leq \alpha \leq l(i)} c_\alpha R_{i\alpha},
c_\alpha\in {\mathbb C}.$ Then
$$
s_i^* vm = s_i^* \gamma(m) v 
= \rho_i (m) s_i^* v
$$
by the intertwining property of $ R_{i\alpha}$, hence
$ s_i^* \gamma(m) vv^* 
= \rho_i (m) s_i^* vv^*
, \forall m\in M$ 
and using the separating property of $\Omega$ for $M$ again we have
$s_i^*\in \Hom(\gamma,\rho_i)$. Since $\rho_i$ are irreducible, 
$\Hom(\rho_i, \rho_j)= \delta_{ij} {\mathbb C}, i,j=1,...,n.$ It follows
that  $\langle s, t\rangle = t^* s$ is an inner product, and
the last part of the lemma follows.
\endproof
Let $s=\sum_{i\beta} s_{i\beta} s_{i\beta}^*$ where $s_{i\beta}$ is the
orthonormal basis as in Lemma \ref{Hil}. When $V$ is zero we set
$s=0$. By construction we have $st=t, \forall t\in V$, in fact
$s$ is the left support of $V$ in $M$.
\begin{lemma}\label{s=1}
$s=1.$
\end{lemma}
\proof
Let us compute
\[
\gamma (v^*) (1-s) v  = \gamma(v^*) v- \sum_{i\beta} \gamma(v^*)
s_{i\beta} s_{i\beta}^* v =
\sum_{i\alpha} m_{i\alpha} R_{i\alpha} - \sum_{i\beta} \gamma(v^*)
s_{i\beta}s_{i\beta}^* v
\]
Note that $\gamma(v^*)s_{i\beta} \in M$,  $s_{i\beta}^* v\in W$ by
Definition \ref{V}. So
\begin{equation}\label{1}
\gamma (v^*) (1-s) v= \sum_{i\alpha} m_{i\alpha}' R_{i\alpha}
\end{equation}
for some $m_{i\alpha}' \in M$. 
If the left hand side above is non-zero, then there is at least
one $m'_{i\alpha}\neq 0$, and we will derive a contradiction by
an averaging trick as on Page 46 of \cite{ILP}. 
Since $M$ is a type III factor, $ m_{i\alpha}' \neq 0$, 
as on Page 46 of \cite{ILP} we can find
$a,b\in M$ such that $E_{\rho_i}(a m_{i\alpha}' \rho_i(b))=1$, where
$E_{\rho_i}: M\rightarrow \rho_i(M)$ is a normal faithful conditional 
expectation. Let $R$ be an injective simple subfactor in $M$ as on 
Page 46 of \cite{ILP} (This  is where the separability of $M_*$ is used), 
and let $u$ be an unitary element of $R$. 
Multiply both sides of   equation (\ref{1}) by $\rho_i(u^*)a$ on the
left and $bu$ on the right, and use the intertwining properties of
$R_{i\alpha},s$ and $v$ we get:
$$
\rho_i(u^*) a \gamma(v^*) \gamma (b) \gamma (u) (1-s) v
= \rho_i(u^*) a m_{i\alpha}' \rho_i(b)  \rho_i(u) R_{i\alpha}
+ \sum_{(j,\beta) \neq (i,\alpha)}  \rho_i(u^*) 
a m_{i\beta}' \rho_j(b)  \rho_j(u) R_{j\beta}.
$$
Averaging by an invariant mean 
over the unitary elements of $R$ (which is amenable since
$R$ is injective) as on Page 46 of \cite{ILP} we get 
$$
t(1-s) v= c_{i\alpha}  R_{i\alpha} + \sum_{(j,\beta) \neq (i,\alpha)}  
c_{j\beta}R_{j\beta},
$$
where $t\in M$, and $c_{j\beta}\in M$ satisfies
$
c_{j\beta} \rho_j(x) = \rho_i(x) c_{j\beta}, \forall x\in R, 
$ 
and $E_{\rho_i} (c_{i\alpha}) = 1$ since $E_{\rho_i}: M\rightarrow 
\rho_i(M)$ is normal. Since $\rho_j, j=1,...,n$ has finite index,
by Cor. 2.11 of \cite{ILP} $c_{j\beta} \in \Hom( \rho_j,  \rho_i)
= \delta_{ij} {\mathbb C}$. Hence we have 
$c_{i\alpha}=1$ and 
$$t(1-s) v= R_{i\alpha} + \sum_{\beta \neq \alpha}
c_{i\beta}R_{i\beta}
$$
where $c_{i\beta} \in {\mathbb C}$. It follows by definition \ref{V}
that $(1-s) t^* \in V$. Since $s$ is the left support of $V$ in $M$, 
we have $s(1-s) t^* = (1-s) t^*$, hence
$(1-s) t^*=0$, and
$ R_{i\alpha} + \sum_{\beta \neq \alpha}
c_{i\beta}R_{i\beta} =0 $
which contradicts the linear independence of $\{R_{j\beta}\}$ over 
${\mathbb C}$. This shows that the left hand side of equation
(\ref{1}) is zero, i.e., 
$\gamma(v^*) (1-s) v =0.$
Multiply this equation on the 
left by $\gamma(m_1)\gamma (v)$ and on the right by $m_2v^*$, and note 
that $(1-s) \in\Hom(\gamma,\gamma)$, we  get
$\gamma (m_1 em_2) (1-s) e=0, \forall m_1,m_2\in M.$ By (1) of
Lemma \ref{basic}, $1$ in the weak closure of $MeM$, therefore
$(1-s)e=0$. Using the fact that $\Omega$ is separating for $M$ and
$e\Omega=\Omega$, we conclude that $s=1$.
\endproof
{\bf The end of the proof of Theorem \ref{key}:}\par
By Lemma \ref{s=1}, $s=\sum_{j\beta} s_{j\beta} s_{j\beta}^*=1,$
and each $ s_{j\beta}\in \Hom(\rho_j,\gamma)$ is an isometry. 
It follows that 
$[\gamma]= \bigoplus_{1\leq j\leq n} l(j)' [\rho_j]$
where $ 0\leq l(j)'< \infty, j=1,...,n.$ By Lemma \ref{can}
$d(\gamma) = \sum_{ 1\leq j \leq n} l(j)' d(\rho_j) < \infty
,$ 
and by Lemma \ref{can} $N\subset M$ has finite index.
\endproof  

\begin{corollary}\label{general} 
Let $N_1,...N_n$ be von Neumann subalgebras of  $M$ where $M$ is a factor 
with separable predual,
$N=N_1\cap N_2...\cap N_n$,
and each $N_i\subset M$ has finite index, $i=1,...,n.$
Assume that:\par
(1) There is a $\phi$ which is a normal faithful state on M
invariant under normal faithful condition expectations $E_i:M\rightarrow N_i$
and $E:M\rightarrow N.$ Let $e_i, e$ be the corresponding Jones projection in
$B(L^2(M,\phi));$\par

(2) $e=e_1\wedge e_2...\wedge e_n.$\par
Then $N\subset M$ has finite index if and only if $e_1,...,e_n$ generate
a finite dimensional algebra.
\end{corollary}
\proof
The only if part follows from Page 2 of \cite{P} and the assumption that
$M$ is a factor. Let us prove the if part. By replacing $M,N_i, e_i,e,\phi$ by
$M\otimes B, N_i\otimes B, e_i\otimes 1, e\otimes 1, \phi\otimes \phi_1$
respectively 
if necessary where $B$ is a type III factor with separable predual and
$\phi_1$ is a normal faithful state on $B$, by (4) of Lemma \ref{basic} it 
is enough to prove the if part by assuming that $M$ is a type III
factor, and $N_i,i=1,2,...,n ,N$ are type III von Neumann algebras.\par

Assume that 
$e_1,...,e_n$ generate
a finite dimensional algebra. By assumption (2) there exists a 
non-commutative polynomial $f$ such that 
$e=f(e_1,...,e_n)$. Let $v_i$ be the isometries in $B(L^2(M,\phi))$
as in Lemma \ref{can} such that $v_i v_i^*= e_i, v_im =\gamma_i(m)v_i
,\forall m\in M.$ Note that  for each $i=1,...,n$ 
$\gamma_i\in \End(M)$ has finite index 
by Lemma \ref{basic} since $N_i\subset M$ has finite index by assumption. 
Using
$ e_i= v_i v_i^*=\gamma_i(v_i^*)v_i, 
\gamma_i(v_i^*) \in M,
$ and the intertwining properties $v_i, i=1,...,n$, it follows that 
$$
e= \sum_{ 1\leq j \leq k} m_j R_j
$$
where $ m_j\in M$ and $R_j\in B(L^2(M,\phi)), 
R_j (m)= \rho_j(m) R_j, \forall m\in M $ and 
$\rho_j\in \End(M)$ is a finite compositions of
$\gamma_i, i=1,...,n$, and hence of finite index by Lemma \ref{dim}.
By Theorem \ref{key} the proof is complete.
\endproof 
\begin{remark}
Note that when $M$ is a type $II_1$ factor in Cor. \ref{general},
one can take $\phi$ to be the trace on $M$, then assumption (1) holds
trivially, and assumption (2) holds by \cite{S}. So we obtain
another proof of Th.\ref{finitecase} when $M$ is a type $II_1$ factor
with separable predual, and this in fact partially inspired  simpler
proofs in section \ref{finite} for finite type von Neumann algebra case.  
We will see a large class of examples in
the case when $M$ is type III where  assumptions (1) and (2) hold
in Cor. \ref{subnet2}.
\end{remark}  
\subsection{Applications to conformal nets}
\subsubsection{Preliminaries on conformal nets}
By an interval of the circle we mean an open connected
non-empty subset $I$ of $S^1$ such that the interior of its 
complement $I'$ is not empty. 
We denote by $\I$ the family of all intervals of $S^1$.

A {\it net} $\A$ of von Neumann algebras on $S^1$ is a map 
\[
I\in\I\to\A(I)\subset B(\H)
\]
from $\I$ to von Neumann algebras on a fixed separable Hilbert space $\H$
that satisfies:
\begin{itemize}
\item[{\bf A.}] {\it Isotony}. If $I_{1}\subset I_{2}$ belong to 
$\I$, then
\begin{equation*}
 \A(I_{1})\subset\A(I_{2}).
\end{equation*}
\end{itemize}
If $E\subset S^1$ is any region, we shall put 
$\A(E)\equiv\bigvee_{E\supset I\in\I}\A(I)$ with $\A(E)=\mathbb C$ 
if $E$ has empty interior (the symbol $\vee$ denotes the von Neumann 
algebra generated). 

The net $\A$ is called {\it local} if it satisfies:
\begin{itemize}
\item[{\bf B.}] {\it Locality}. If $I_{1},I_{2}\in\I$ and $I_1\cap 
I_2=\emptyset$ then 
\begin{equation*}
 [\A(I_{1}),\A(I_{2})]=\{0\},
 \end{equation*}
where brackets denote the commutator.
\end{itemize}
The net $\A$ is called {\it M\"{o}bius covariant} if in addition 
satisfies
the following properties {\bf C,D,E,F}:
\begin{itemize}
\item[{\bf C.}] {\it M\"{o}bius covariance}. 
There exists a non-trivial strongly 
continuous unitary representation $U$ of the M\"{o}bius group 
$\Mob$ (isomorphic to $PSU(1,1)$) on $\H$ such that
\begin{equation*}
 U(g)\A(I) U(g)^*\ =\ \A(gI),\quad g\in \Mob,\ I\in\I.
\end{equation*}
\item[{\bf D.}] {\it Positivity of the energy}.  
The generator of the one-parameter
rotation subgroup of $U$ (conformal Hamiltonian), denoted by 
$L_0$ in the following,  is positive.
\item[{\bf E.}] {\it Existence of the vacuum}.  There exists a unit
$U$-invariant vector $\Omega\in\H$ (vacuum vector), and $\Omega$ is
cyclic for the von Neumann algebra $\bigvee_{I\in\I}\A(I)$.
\end{itemize}
By the Reeh-Schlieder theorem $\Omega$ is cyclic and separating for 
every fixed $\A(I)$. The modular objects associated with 
$(\A(I),\Omega)$ have a geometric meaning
\[
\Delta^{it}_I = U(\Lambda_I(2\pi t)),\qquad J_I = U(r_I)\ .
\]
Here $\Lambda_I$ is a canonical one-parameter subgroup of $\Mob$ and $U(r_I)$ is a 
antiunitary acting geometrically on $\A$ as a reflection $r_I$ on $S^1$. 

This implies {\em Haag duality}: 
\[
\A(I)'=\A(I'),\quad I\in\I\ ,
\]
where $I'$ is the interior of $S^1\setminus I$.

\begin{itemize}
\item[{\bf F.}] {\it Irreducibility}. $\bigvee_{I\in\I}\A(I)=B(\H)$. 
Indeed $\A$ is irreducible iff
$\Om$ is the unique $U$-invariant vector (up to scalar multiples). 
Also  $\A$ is irreducible
iff the local von Neumann 
algebras $\A(I)$ are factors. In this case they are III$_1$-factors 
with separable predual
in 
Connes classification of type III factors.
\end{itemize}
By a {\it conformal net} (or diffeomorphism covariant net)  
$\A$ we shall mean a M\"{o}bius covariant net such that the following 
holds:
\begin{itemize}
\item[{\bf G.}] {\it Conformal covariance}. There exists a projective 
unitary representation $U$ of $\Diff(S^1)$ on $\H$ extending the unitary 
representation of $\Mob$ such that for all $I\in\I$ we have
\begin{gather*}
 U(\phi)\A(I) U(\phi)^*\ =\ \A(\phi.I),\quad  \phi\in\Diff(S^1), \\
 U(\phi)xU(\phi)^*\ =\ x,\quad x\in\A(I),\ \phi\in\Diff(I'),
\end{gather*}
\end{itemize}
where $\Diff(S^1)$ denotes the group of smooth, positively oriented 
diffeomorphism of $S^1$ and $\Diff(I)$ the subgroup of 
diffeomorphisms $g$ such that $\phi(z)=z$ for all $z\in I'$.
\par
Note that if $\phi\in \Diff(I)$, then $U(\phi)\in \A(I)$ by {\bf G}
and Haag duality.\par
Next we  recall some definitions from \cite{KLM} .
Recall that   ${\I}$ denotes the set of intervals of $S^1$.
Let $I_1, I_2\in {\I}$. We say that $I_1, I_2$ are disjoint if
$\bar I_1\cap \bar I_2=\emptyset$, where $\bar I$
is the closure of $I$ in $S^1$.  
When $I_1, I_2$ are disjoint, $I_1\cup I_2$
is called a 1-disconnected interval in \cite{Xu3}.  
Denote by ${\I}_2$ the set of unions of disjoint 2 elements
in ${\I}$. Let ${\A}$ be an irreducible M\"{o}bius covariant net
. For $E=I_1\cup I_2\in{\I}_2$, let
$I_3\cup I_4$ be the interior of the complement of $I_1\cup I_2$ in 
$S^1$ where $I_3, I_4$ are disjoint intervals. 
Let 
$$
{\A}(E):= A(I_1)\vee A(I_2), \quad
\hat {\A}(E):= (A(I_3)\vee A(I_4))'.
$$ Note that ${\A}(E) \subset \hat {\A}(E)$.
Recall that a net ${\A}$ is {\it split} if ${\A}(I_1)\vee
{\A}(I_2)$ is naturally isomorphic to the tensor product of
von Neumann algebras ${\A}(I_1)\otimes
{\A}(I_2)$ for any disjoint intervals $I_1, I_2\in {\I}$.
${\A}$ is {\it strongly additive} if ${\A}(I_1)\vee
{\A}(I_2)= {\A}(I)$ where $I_1\cup I_2$ is obtained
by removing an interior point from $I$.
\bdefin\label{rational}
\cite{{KLM}, {LX}}
${\A}$ is said to be completely  rational if
${\A}$ is split, and 
the index $[\hat {\A}(E): {\A}(E)]$ is finite for some
$E\in {\I}_2$ . The value of the index
$[\hat {\A}(E): {\A}(E)]$ (it is independent of 
$E$ by Prop. 5 of \cite{KLM}) is denoted by $\mu_{{\A}}$
and is called the $\mu$-index of ${\A}$. If 
the index $[\hat {\A}(E): {\A}(E)]$ is infinity for some
$E\in {\I}_2$, we define the $\mu$-index of ${\A}$ to be
infinity.
%\label{Definition 2.2}
\edefin
Note that, by recent results in \cite{LX}, 
every irreducible, split, 
local conformal net with finite $\mu$-index is automatically strongly 
additive. Hence we have modified the definition in \cite{KLM} by dropping the
strong additivity requirement in the above definition. Also note that
if $\A$ is completely rational, then $\A$ has only finitely many 
irreducible covariant representations by \cite{KLM}.
\par
Let $\A$ be a conformal net. By a {\it conformal subnet} (cf. \cite{L2}) 
we shall mean a
map  
\[
I\in\I\to\B(I)\subset \A(I)
\]
that associates to each interval $I\in \I$ a von Neumann subalgebra $\B(I)$
of $\A(I)$, which is isotonic
\[
\B(I_1)\subset \B(I_2), I_1\subset I_2,
\]
and   M\"{o}bius covariant with respect to the the representation $U$, 
namely 
\[
U(g) \B(I) U(g)^*= \B(g.I)
\] for all $g\in \Mob$ and $I\in \I$. Note that by Lemma 13 
of \cite{L2} for each $I\in \I$ there exists a conditional 
expectation $E_I: \A(I)\rightarrow \B(I)$ such that $E$ preserves the
vector state given by the vacuum of $\A$.

\begin{lemma}\label{irred1}
Let $\B\subset \A$ be a conformal subnet and assume that
$U(\Diff(I))\subset \B(I), \forall I\in \I$, and $\B$ is 
completely rational.  Then: \par
(1) $\B\subset \A$ is irreducible, i.e., 
$\B(I)'\cap \A(I) = {\mathbb C},  \forall I\in \I$;\par
(2) $\B(I)\subset \A(I)$
has finite index $\forall I\in \I$, and $\A$ is 
completely rational.
\end{lemma}
\proof
Let $p\in \B(I)'\cap \A(I)$. Since $\B(I')\subset \A(I')$, we
have $p\in (\B(I)\vee \B(I'))'$. Since $\B$ is completely rational, 
by strong additivity $p\in (\vee_{I\in \I} \B(I))'$, and so
$p\in (\vee_{I\in \I} U(\Diff(I)))'$ by assumption. Notice that 
$\Diff(I), I\in \I$ generates $\Diff (S^1)$, and so 
$p$ commutes with $  U(\Diff (S^1))$, and so $p\Omega$ is an eigenvector
of the conformal Hamiltonian with eigenvalue $0$. By {\bf F} it follows
that $p\Omega = x\Omega, x\in {\mathbb C}$, therefore $p= x$ since
$\Omega$ is separating for $\A(I)$, proving (1). (2) follows from 
(1) , Prop. 2.3 of \cite{KL} and Th. 24 of \cite{L2}.
\endproof
\begin{corollary}\label{subnet1}
Let $\A$ be a conformal net, and $\B_1, \B_2,...\B_n$ be conformal
subnets of $\A$. Let $\B$ be the conformal
subnet of $\A$ such that $\B(I)= \B_1(I)\cap  \B_2(I)...\cap \B_n(I), 
\forall I\in \I$. Assume that: \par
(1) $U(\Diff(I))\subset \B_i(I), \forall I\in \I, i=1,...,n$; \par
(2) Each $\B_i$ is completely rational,  $i=1,...,n$; \par
(3) $e=e_1\wedge e_2\cdots\wedge e_n$, where 
$e$ (resp. $e_i, i=1,...,n$) are Jones projections 
onto $\overline{\B(I)\Omega}$
(resp.  $\overline{\B_i(I)\Omega},i=1,...,n $).\par
Then $\B$ is completely rational if and only if 
$e_1,...,e_n$ generate a finite dimensional algebra.
\end{corollary}
\proof
Note that by definition and assumption (1) we have
$U(\Diff(I))\subset \B(I), \forall I\in \I$. By Lemma
\ref{irred1} and assumptions (1) and (2) each 
$\B_i(I) \subset \A(I)$ has finite index, and $\A$ is
completely rational. By Lemma \ref{irred1} again we have that
 $\B$ is completely rational if and only if $\B(I)\subset \A(I)$
has finite index. We note that $\A(I)$ is a type $III_1$ factor 
with separable predual as
stated in {\bf F}. The assumption (1) of Corollary \ref{general}
in this case follows by the remark after the definition of
conformal subnets. 
The corollary now follows from
Corollary \ref{general}.
\endproof
A large class of conformal subnets verifying assumptions (1), (2) and
(3) of Cor. \ref{subnet1} come from cosets and orbifolds (cf. \cite{Xu1}
and \cite{Xu2}). Let us recall some definitions. 
Let $G$ be a simply connected  compact Lie group. By Th. 3.2
of \cite{FG}, 
the vacuum positive energy representation of the loop group
$LG$ (cf. \cite{PS}) at level $k$ on a Hilbert space $\H$, denoted by
$\pi^0$,  
gives rise to an irreducible conformal net 
denoted by {\it ${\A}_{G}$} when $k$ is fixed. We will use
$\Omega$ to denote the vacuum vector. Note that
$ {\A}_{G}(I) = \pi^0(L_IG)''$ where $L_IG$ are these elements of
$LG$ which are equal to identity of $G$ on $I'$. 
By Th. 13.4.2 of \cite{PS} (also cf. \cite{GW}) there is a projective 
unitary representation $U$ of $\Diff(S^1)$ on $\H$ such that
$U(\phi) \pi^0(f) U(\phi)^*= \pi^0(f\cdot \phi^{-1})$ for any
$f\in LG$.\par
Let $H\subset G$ be a simply connected Lie subgroup. We define 
a conformal subnet $\A_{H,G/H}$ of $\A$ by 
$\A_{H,G/H}(I):= \pi^0(L_IH)''\vee  (\pi^0(L_IH)'\cap  \pi^0(L_IG)'')$.
Recall from \S3 of \cite{Xu1}  that $H\subset G$ is {\it cofinite}
if $\A_{H,G/H}(I)\subset {\A}_{G}(I)$ has finite index for some $I\in \I$.
Note that Conjecture 2.13 of \cite{Xu1} implies that any such 
$H\subset G$ is cofinite. See Cor. 3.4 of \cite{Xu1} for a list of
inclusions which have been proved to be cofinite.\par
When $\Gamma$ is a finite subgroup of $G$, we denote by 
$\A^\Gamma$ a conformal subnet of ${\A}_{G}$ such that
$\A^\Gamma (I)= \{ m\in {\A}_{G}(I) | \pi^0(h)m = m \pi^0(h), \forall 
h\in \Gamma\}$. This is an example of orbifold construction in 
\cite{Xu2}. \par
\begin{lemma}\label{irred2}
Let $H_1,..., H_l$ be simply connected Lie 
subgroups of $G$, and let $\Gamma_{l+1},...\Gamma_{n}$ be finite
subgroups of $G$.
Let $\B_i=\A_{H_i,G/H_i}, i=1,...,l$ and
$\B_j= \A^{\Gamma_j}, j=l+1,...,n$ be conformal subnets of $ {\A}_{G}$
as above, and  let $\B$ be the conformal
subnet of $\A$ such that $\B(I)= \B_1(I)\cap  \B_2(I)...\cap \B_n(I), 
\forall I\in \I$. Then:\par
(1) $U(\Diff(I))\subset \B_i(I), \forall I\in \I, i=1,...,n$;\par
(2) $e=e_1\wedge e_2\cdots\wedge e_n$, where 
$e$ (resp. $e_i, i=1,...,n$) are Jones projections 
onto $\overline{\B(I)\Omega}$
(resp.  $\overline{\B_i(I)\Omega},i=1,...,n $). 
\end{lemma}
\proof
Ad (1):
When $\B_j= \A^{\Gamma_j}, j=l+1,...,n$ (1) holds trivially since
$U(\Diff (S^1))$ commutes with $\Gamma_j$. Let us assume that
$\B_i=\A_{H_i,G/H_i}$. Let $\phi\in \Diff (I)$. By Th. 13.4.2 
of \cite{PS} (also cf. \cite{GW}) we have 
$U(\phi) \pi^0(f) U(\phi)^*= \pi^0(f\cdot \phi^{-1})$ for any
$f\in LG$, and by the remark after {\bf G} $U(\phi)\in \A_G(I)$. Apply
the same argument to $LH$ we conclude that there is a unitary 
element $\tilde\phi
\in \pi^0(L_IH)''$ such that 
$\tilde\phi \pi^0(f) \tilde\phi^*= 
\pi^0(f\cdot \phi^{-1})
= U(\phi)\pi^0(f) U(\phi)^*
$ for any
$f\in L_IH$. It follows that 
$\tilde\phi^* U(\phi) \in  \pi^0(L_IH)'\cap  \pi^0(L_IG)''$,
hence $ U(\phi)\in \A_{H_i,G/H_i}(I)= \B_i$.\par
Ad (2):
Note by definition 
$\overline{B(I)\Omega}\subset \cap_{1\leq i\leq n} \overline{B(I)\Omega}
$, and so it is sufficient to show that for any eigenvector $\psi$ of
$L_0$ with eigenvalue $m\geq 0$ in 
$\cap_{1\leq i\leq n} \overline{B(I)\Omega}$,  $\psi$ is also in 
$\overline{B(I)\Omega}$. By Reeh-Schlieder theorem it is sufficient to 
show that $\psi\in  \overline{\vee_{I\in \I}B(I)\Omega}$. The proof
is essentially contained on Page 22 of \cite{Xu1} as follows:
Choose
two smooth functions $f_1(z)$ and $ f_2(z)$ on the 
unit circle,  with support $f_1\subset I_1
\in \I$, support $f_2\subset I_2
\in \I$ and $f_1+f_2= 1$. Then 
$\psi = V(\psi, z^{-1})\Omega =V(\psi, z^{-1}f_1)\Omega+ 
V(\psi, z^{-1}f_2)\Omega  
$
where $ V(\psi,\cdot)$ are the smeared vertex operators as
defined on Page 11 of \cite{Xu1}. By the same proof as in Prop. 2.11 of
\cite{Xu1}, $V(\psi, z^{-1}f_1)$ is a closed operator affiliated with 
with $\B_i(I_1), i=1,2,...,n$, and so $V(\psi, z^{-1}f_1)$  
is a closed operator affiliated with  $\B(I_1)$, it follows that 
$V(\psi, z^{-1}f_1)\Omega\in \overline{B(I_1)\Omega}$. Similarly
$V(\psi, z^{-1}f_2)\Omega\in \overline{B(I_2)\Omega}$, and 
we conclude that  $\psi\in  \overline{\vee_{I\in \I}B(I)\Omega}$.
\endproof
\begin{remark}
Due to Lemma \ref{irred2}, we conjecture that assumption (3) of
Cor. \ref{subnet1} is always satisfied.
\end{remark}
We note that by \cite{Xu3} which is based on 
\cite{W}, $\A_G$ is completely rational 
if $G=SU(N_1)\times SU(N_2)\times...\times SU(N_k)$, and it has been 
conjectured 
that
$\A_G$ is completely rational for all $G$.\par
\begin{corollary}\label{subnet2}
Let $H_1,..., H_l$ be simply connected Lie 
subgroups of $G$, and let $\Gamma_{l+1},...\Gamma_{n}$ be finite
subgroups of $G$.
Let $\B_i=\A_{H_i,G/H_i}, i=1,...,l$ and
$\B_j= \A^{\Gamma_j}, j=l+1,...,n$ be conformal subnets of $ {\A}_{G}$
as described before Lemma \ref{irred2}, and 
 let $\B$ be the conformal
subnet of $\A$ such that $\B(I)= \B_1(I)\cap  \B_2(I)...\cap \B_n(I), 
\forall I\in \I$. Assume that  each $H_i \subset G, i=1,2,...,l$ is cofinite
and  $\A_G$ is completely rational.
Then $\B(I)$ is completely rational if and only if 
$e_1,...,e_n$ generate a finite dimensional algebra where
$e_i, i=1,...,n$ are Jones projections from $\H$ 
onto $\overline{\B_i(I)\Omega},i=1,...,n .$
\end{corollary}
\proof
Note that since we assume that  each $H_i \subset G, i=1,2,...,l$ is cofinite,
$\B_i(I)\subset \A_G(I)$ has finite index for $i=1,...,n$. By 
Th. 24 of \cite{L2} each $\B_i, i=1,...,n$ is completely rational. Hence
assumption (2) of Cor. \ref{subnet1} is satisfied. Note that assumptions
(1) and (3) of Cor. \ref{subnet1} 
are satisfied thanks to Lemma \ref{irred2}. Hence the 
corollary is proved by Cor.\ref{subnet1}.
\endproof
We note that the nature of the algebra generated by 
$e_1,...,e_n$ in Cor. \ref{subnet2} can in principle be determined by
the representation theory information about pairs $LH_i\subset LG$
and $\Gamma_j\subset LG$. In the case $n=2$, it is well-known that
$e_1,e_2$ generate a finite dimensional algebra if and only if
the ``angle operator'' $e_1e_2e_1$ has finite spectrum. It is
an interesting question to see if one can obtain new examples of
completely rational conformal nets by using Cor. \ref{subnet2}.    
\end{document}